\documentclass{birkjour}

 \newtheorem{thm}{Theorem}[section]

 \newtheorem{prop}[thm]{Proposition}
 \theoremstyle{definition}
 
 \theoremstyle{remark}

\begin{document}
%
%
%
%
%
%
%
%
%
\title[Hypercomplex structures with Hermitian-Norden metrics]
      {Hypercomplex structures with Hermitian-\\Norden metrics
    on four-dimensional Lie \\algebras}

\author[M. Manev]{Mancho Manev}

\address{%
Paisii Hilendarski University of Plovdiv \br Faculty of Mathematics, Informatics and IT\br
Department of Algebra and Geometry \br 236 Bulgaria blvd \br
Plovdiv 4027 \br Bulgaria}

\email{mmanev@uni-plovdiv.bg}

\thanks{This paper is partially supported by project NI13-FMI-002
of the Scientific Research Fund, Plovdiv University, Bulgaria
and the German Academic Exchange Service (DAAD)}

\subjclass{53C15, 53C50, 22E30, 53C55, 53C56}

\keywords{hypercomplex structure, 4-dimensional Lie algebra, Hermitian metric, Norden metric,
indefinite metric}


\begin{abstract}
Integrable hypercomplex structures with Hermitian and Norden metrics on Lie groups of dimension 4
are considered. The corresponding five types of invariant hypercomplex structures with
hyper-Hermitian metric, studied by M.L. Barberis, are constructed here. The different cases
regarding the signature of the basic pseudo-Riemannian metric are considered.
\end{abstract}


\newcommand{\X}{\mathfrak{X}}
\newcommand{\K}{\mathcal{K}}
\newcommand{\s}{\mathfrak{S}}
\newcommand{\g}{\mathfrak{g}}
\newcommand{\z}{\mathfrak{z}}
\newcommand{\R}{\mathbb{R}}
\newcommand{\C}{\mathbb{C}}
\newcommand{\W}{\mathcal{W}}
\newcommand{\HC}{\mathcal{H}}

\newcommand{\n}{\nabla}

\newcommand{\ie}{i.e. }
\newcommand{\norm}[1]{\left\Vert#1\right\Vert ^2}
\newcommand{\nN}{\norm{N}}
\newcommand{\tr}{{\rm tr}}

\newcommand{\nJ}[1]{\norm{\nabla J_{#1}}}
\newcommand{\thmref}[1]{Theorem~\ref{#1}}
\newcommand{\propref}[1]{Proposition~\ref{#1}}
\newcommand{\secref}[1]{\S\ref{#1}}
\newcommand{\lemref}[1]{Lemma~\ref{#1}}
\newcommand{\corref}[1]{Corollary~\ref{#1}}
\newcommand{\dfnref}[1]{Definition~\ref{#1}}


\newcommand{\al}{\alpha}
\newcommand{\bt}{\beta}
\newcommand{\gm}{\gamma}
\newcommand{\om}{\omega}
\newcommand{\lm}{\lambda}
\newcommand{\ta}{\theta}
\newcommand{\ea}{\varepsilon_\alpha}
\newcommand{\eb}{\varepsilon_\bt}
\newcommand{\eg}{\varepsilon_\gamma}

\maketitle

\section*{Introduction}

The present work is inspired by the work of Barberis \cite{Barb}
where invariant hypercomplex structures $H$ on 4-dimensional
real Lie groups are classified.
In that case the corresponding metric is positive definite and
Hermitian with respect to the triplet of complex structures of $H$.
Our main goal is to classify 4-dimensional real Lie algebras
which admit hypercomplex structures with Hermitian and Norden metrics.

We equip a hypercomplex structure $H$ with a metric structure,
generated by a pseudo-Riemannian metric $g$ of neutral signature
(see \cite{GrMa,GrMaDi}). In our case the one (resp., the other
two) of the almost complex structures of $H$ acts as an isometry
(resp., act as anti-isometries) with respect to $g$ in each
tangent fibre.
Thus,  there exist three (0,2)-tensors
associated by $H$ except the metric $g$
--- a K\"ahler form and two metrics of the same type.
The metric $g$ is Hermitian with respect to the one almost
complex structure of $H$ and $g$ is a Norden (\ie an anti-Hermitian)
metric regarding the other two almost complex structures
of $H$. For this reason we call the derived almost hypercomplex
structure an \emph{almost hypercomplex structure with Hermitian-Norden-metrics}
or briefly \emph{almost hypercomplex HN-metric structure}.

Let us remark that in \cite{Snow} and \cite{Ovan} are classified
the invariant complex structures on 4-dimensional solvable
simply-connected real Lie groups where the dimension of commutators
is less than three and equal three, respectively.

A hypercomplex structure is called \emph{Abelian} (\cite{BaDoMi}) if
$[J_{\al}x, J_{\al}y ] = [x, y]$,
for all $x, y \in \g$ ($\al = 1, 2, 3$). Abelian hypercomplex structures are
considered in \cite{BaDo}, \cite{DotFin} and
they can only occur on solvable Lie algebras (\cite{FinGran}).
It is clear that the condition %
$N(x,y) = [Jx,Jy]-J[Jx,y]-J[x,Jy]-[x,y]=0$ can be rewritten as
$[Jx,Jy]-[x,y] = J([Jx,y]+[x,Jy])$ for all $x, y \in \g$. Thus,
Abelian complex structures and therefore Abelian hypercomplex
structure are integrable.

If the three almost complex structures of $H$ are parallel with
respect to the Levi-Civita connection $\nabla$ of $g$ then such
hypercomplex HN-manifolds of K\"ahler type we call
\emph{hyper-K\"ahler HN-manifolds}, which are flat according to
\cite{GrMaDi}.

The paper is organized as follows. In Sect.~\ref{sec1} we recall some facts
about the almost hypercomplex HN-manifolds known from
\cite{AlMa,GrMa,GrMaDi,Ma09}.
In Sect.~\ref{sec2} we construct different types of hypercomplex
structures on Lie algebras following the Barberis classification.
%

The basic problem of this work is the existence and the geometric
characteristics of hypercomplex HN-structures on 4-dimensional
Lie algebras according to the Barberis classification. The
main results of this paper is construction of the different types
of the considered structures and their characterization.

\section{Preliminaries}\label{sec1}

Let $(M,H)$ be a hypercomplex manifold, \ie $M$ is a
$4n$-dimension\-al differentiable manifold and $H=(J_1,J_2,J_3)$
is a triple of complex structures on $M$ with the following
properties for all cyclic permutations $(\al, \bt, \gm)$ of
$(1,2,3)$:%
\begin{equation}\label{J123} %
J_\al=J_\bt\circ J_\gm=-J_\gm\circ J_\bt, \qquad
J_\al^2=-I,
\end{equation} %
where $I$ denotes the identity and 
moreover, it is valid %
\begin{equation}\label{N=0} %
N_\al=0,\qquad \al\in\{1,2,3\}
\end{equation} %
for the Nijenhuis tensors $N_\al$ of $J_\al$ given by
\begin{equation}\label{N_al} %
N_\al(\cdot,\cdot)\allowbreak{}= \left[J_\al \cdot,J_\al \cdot
\right]
    -J_\al\left[J_\al \cdot,\cdot \right]
    -J_\al\left[\cdot,J_\al \cdot \right]
    -\left[\cdot,\cdot \right]
\end{equation}%
on $\X(M)$.

A hypercomplex structure on a Lie group $G$ is said to be \emph{invariant}
if left translations by elements of $G$ are holomorphic with respect to
$J_{\al}$ for all $\al\in\{1,2,3\}$. Obviously, if $\g$ is
the corresponding Lie algebra of the Lie group $G$, a
hypercomplex structure on $\g$ induces an invariant hypercomplex structure
on $G$ by left translations.

Let $G$ be a simply connected 4-dimensional real Lie group
admitting an invariant hypercomplex structure. A left invariant
metric on $G$ is called invariant hyper-Hermitian  if it is
hyper-Hermitian with respect to some invariant hypercomplex
structure on $G$. It is known that all such metrics on given $G$
are equivalent up to homotheties.

If $\g$ denotes the Lie algebra of $G$ then it is known the following

\begin{thm}[\cite{Barb}]\label{thm-Barb}
The only 4-dimensional Lie algebras admitting an integrable
hypercomplex structure are the following types:

\emph{\textbf{(hc1)}} $\g$ is Abelian;$\qquad\qquad$
\emph{\textbf{(hc2)}} $\g\cong\R\oplus\mathfrak{so}(3)$;
$\qquad\qquad$
\emph{\textbf{(hc3)}} $\g\cong\mathfrak{aff}(\C)$;

\emph{\textbf{(hc4)}} $\g$ is the solvable Lie algebra corresponding to $\R H^4$;

\emph{\textbf{(hc5)}} $\g$ is the solvable Lie algebra corresponding to $\C H^2$,
\\
where $\R\oplus\mathfrak{so}(3)$ is the Lie algebra of the Lie
groups $U(2)$ and $S^3\times S^1$; $\mathfrak{aff}(\C)$ is the Lie
algebra of the affine motion group of $\C$ -- the unique
4-dimensional Lie algebra carrying an Abelian hypercomplex
structure; $\R H^4$ is the real hyperbolic space; $\C H^2$ is the
complex hyperbolic space.
\end{thm}

Let $g$ be a neutral metric on $(M,H)$ with the properties
\begin{equation}\label{gJJ} %
g(\cdot,\cdot)=\ea g(J_\al \cdot,J_\al \cdot), \end{equation} %
where
\[
\ea=
\begin{cases}
\begin{array}{ll}
1, \quad & \al=1;\\
-1, \quad & \al=2;3.
\end{array}
\end{cases}
\]
Moreover, the associated (K\"ahler) 2-form $g_1$ and the
associated neutral metrics $g_2$ and $g_3$ are determined by
\begin{equation}\label{gJ} %
g_\al(\cdot,\cdot)=g(J_\al \cdot,\cdot)=-\ea g(\cdot,J_\al \cdot).
\end{equation}%

The structure tensors of a such manifold are the following three
$(0,3)$-tensors
\begin{equation}\label{F}
F_\al (x,y,z)=g\bigl( \left( \n_x J_\al
\right)y,z\bigr)=\bigl(\n_x g_\al\bigr) \left( y,z \right),
\end{equation}
where $\n$ is the Levi-Civita connection generated by $g$.
The corresponding Lee 1-forms $\ta_\al$ are defined by
\begin{equation}\label{theta-al}
\ta_\al(\cdot)=g^{ij}F_\al(e_i,e_j,\cdot)
\end{equation}%
for an arbitrary basis $\{e_1,e_2,\dots, e_{4n}\}$ of $T_pM$,
$p\in M$.

In \cite{GrMaDi} we study the so-called \emph{hyper-K\"ahler
manifolds with HN-metric structure} (or pseudo-hyper-K\"ahler
manifolds), \ie the almost hypercomplex HN-manifold in the class
$\K$, where $\n J_\al=0$ for all $\al=1,2,3$. A sufficient
condition $(M,H,G)$ be in $\K$ is this manifold be of
K\"ahler-type with respect to two of the three complex structures
of $H$ \cite{GrMa}.
%

As $g$ is an indefinite metric, there exist isotropic vectors $x$
on $M$, \ie{} \(g(x,x)=0\), \(x\neq 0\). 
In \cite{GrMa} we define the invariant square norm
\begin{equation}\label{nJ}
\nJ{\al}= g^{ij}g^{kl}g\bigl( \left( \nabla_i J_\al \right) e_k,
\left( \nabla_j J_\al \right) e_l \bigr),
\end{equation}
where $\{e_1,e_2,\dots, e_{4n}\}$ is an arbitrary basis of the
tangent space $T_pM$ at an arbitrary point $p\in M$ of $T_pM$. We
say that an almost hypercomplex HN-manifold is an \emph{isotropic
hyper-K\"ahler HN-manifold} if $\nJ{\al}=0$ for each $J_\al$ of
$H$. Clearly, if the manifold is a hyper-K\"ahler
HN-manifold, then it is an isotropic hyper-K\"ahler HN-manifold.
The inverse statement does not hold.

Let us note that according to \eqref{gJJ} the manifold $(M,J_1,g)$
is almost Hermitian and the manifolds $(M,J_\al,g)$, $\al=2,3$,
are almost complex manifolds with Norden metric \cite{GaBo}. The
basic classes of the mentioned two types of manifolds are given in
\cite{GrHe} and \cite{GaBo}, respectively. The special class
$\W_0(J_\al):$ $F_\al=0$ $(\al =1,2,3)$ of the K\"ahler-type
manifolds belongs to any other class within the corresponding
classification. In the 4-dimensional case the four basic classes
of the almost Hermitian manifolds are restricted to two:
$\W_2(J_1)$, the class of the almost K\"ahler manifolds and
$\W_4(J_1)$, the class of the Hermitian manifolds with respect to
$J_1$. They are determined for $\dim M=4$ by:
\begin{equation}\label{cl-H}
\begin{split}
&\W_2(J_1):\; \mathop{\s}_{x,y,z}\bigl\{F_1(x,y,z)\bigr\}=0; \\
&\W_4(J_1):\; F_1(x,y,z)=\frac{1}{2}
                \left\{g(x,y)\ta_1(z)-g(x,z)\ta_1(y)\right. \\
&\phantom{\W_4(J_1):\; F_1(x,y,z)=\frac{1}{2}
                \left\{\right.}
                \left.-g(x,J_1y)\ta_1(J_1z)+g(x,J_1z)\ta_1(J_1y)
                \right\},
\end{split}
\end{equation}
where $\s $ is the cyclic sum by three
arguments $x$, $y$, $z$.
The basic classes of the almost Norden manifolds (i.e., for $\al=2$ or $3$)
are determined for dimension $4$ as follows:
\begin{equation}\label{cl-N}
\begin{split}
&\W_1(J_\al):\; F_\al(x,y,z)=\frac{1}{4}\bigl\{
g(x,y)\ta_\al(z)+g(x,z)\ta_\al(y)\bigr.\\
&\phantom{\W_1(J_\al):\; F_\al(x,y,z)=\frac{1}{4}\bigl\{\bigr.} %
\bigl.+g(x,J_\al y)\ta_\al(J_\al z)
    +g(x,J_\al z)\ta_\al(J_\al y)\bigr\};\\
&\W_2(J_\al):\; \mathop{\s}_{x,y,z}
\bigl\{F_\al(x,y,J_\al z)\bigr\}=0,\qquad \ta_\al=0;\\
&\W_3(J_\al):\; \mathop{\s}_{x,y,z} \bigl\{F_\al(x,y,z)\bigr\}=0.
\end{split}
\end{equation}

It is known that the class of the complex manifolds with Norden metric is
$\W_1\oplus\W_2$ for $J_{\al}$ ($\al=2,3$).

Then the class of hypercomplex manifolds with Hermitian-Norden metrics is
\[
\HC=\W_4(J_1)\cap\left(\W_1\oplus\W_2\right)(J_2)\cap\left(\W_1\oplus\W_2\right)(J_3).
\]

\section{Four-dimensional Lie algebras with such structures}\label{sec2}

Let $\{e_1,e_2,e_3,e_4\}$ be a basis of a 4-dimensional real Lie
algebra $\g$ with center $\z$ and derived Lie algebra
$\g'=[\g,\g]$. A standard hypercomplex structure on $\g$ is
defined as in \cite{So}:
\begin{equation}\label{JJJ}
\begin{array}{llll}
J_1e_1=e_2, \quad & J_1e_2=-e_1,  \quad &J_1e_3=-e_4, \quad &J_1e_4=e_3;
\\[4pt]
J_2e_1=e_3, &J_2e_2=e_4, &J_2e_3=-e_1, &J_2e_4=-e_2;
\\[4pt]
J_3e_1=-e_4, &J_3e_2=e_3, &J_3e_3=-e_2, &J_3e_4=e_1.\\[4pt]
\end{array}
\end{equation}

Let us introduced a pseudo-Euclidian metric $g$ with
neutral signature as follows
\begin{equation}\label{g}
g(x,y)=x^1y^1+x^2y^2-x^3y^3-x^4y^4,
\end{equation}
where $x(x^1,x^2,x^3,x^4)$, $y(y^1,y^2,y^3,y^4) \in \g$. This
metric satisfies \eqref{gJJ} and \eqref{gJ}. Then the metric $g$
generates an almost hypercomplex HN-metric structure on $\g$.

Let us consider the different cases of \thmref{thm-Barb}.

\subsection{Hypercomplex HN-metric structure of  type (hc1)}
Obviously, in this case the considered manifold belongs to the class $\K$.

\subsection{Hypercomplex HN-metric structure of type (hc2)}

Let $\g$ be not solvable and let us determine it by %
\begin{equation}\label{HC2a}
[e_2,e_4]=e_3,\qquad
[e_4,e_3]=e_2,\qquad [e_3,e_2]=e_4. %
\end{equation}
In this consideration the $(+)$-unit $e_1\in\R$, \ie $g(e_1,e_1)=1$,
is orthogonal to $\g'$ with respect to $g$.

Then we compute covariant derivatives in the basis and the nontrivial ones are
\begin{equation}\label{HC2a-n}
\begin{array}{lll}
\n_{e_2} e_3=-\frac{3}{2}e_4,\qquad & \n_{e_3}
e_2=-\frac{1}{2}e_4,\qquad &
\n_{e_4} e_2=\frac{1}{2}e_3, \\[4pt]
\n_{e_2} e_4=\frac{3}{2}e_3,\qquad & \n_{e_3}
e_4=-\frac{1}{2}e_2,\qquad & \n_{e_4} e_3=\frac{1}{2}e_2.
\end{array}
\end{equation}

By virtue of \eqref{HC2a-n}, \eqref{JJJ} and  \eqref{F}, we obtain components
$(F_{\al})_{ijk}=F_{\al}(e_i,e_j,e_k)$, $i,j,k\in\{1,2,3,4\}$, as follows:
\begin{equation}\label{HC2aFa}
\begin{aligned}
&
\begin{aligned}
(F_{1})_{314}&=-(F_{1})_{323}=(F_{1})_{332}=-(F_{1})_{341}=
\\[4pt]
&=-(F_{1})_{413}=-(F_{1})_{424}=(F_{1})_{431}=(F_{1})_{442}=\frac{1}{2};\\[4pt]
\end{aligned}\\ &
\begin{aligned}
(F_{2})_{214}&=-(F_{2})_{223}=-(F_{2})_{232}=(F_{2})_{241}=\frac{3}{2},\quad
&(F_{2})_{322}=-1,
\\[4pt]
(F_{2})_{412}&=(F_{2})_{421}=(F_{2})_{434}=(F_{2})_{443}=\frac{1}{2}, \quad &(F_{2})_{344}=-1;\\[4pt]
(F_{3})_{213}&=(F_{3})_{224}=(F_{3})_{231}=(F_{3})_{242}=\frac{3}{2},\quad
&(F_{3})_{422}=1,
\\[4pt]
(F_{3})_{312}&=(F_{3})_{321}=-(F_{3})_{334}=-(F_{3})_{343}=\frac{1}{2},
\quad &(F_{3})_{433}=1.
\end{aligned}
\end{aligned}
\end{equation}
 The only non-zero components $(\theta_{\al})_i=(\theta_{\al})(e_i)$, $i=1,2,3,4$,
 of the corresponding Lee forms are
\begin{equation}\label{HC2a-titi}
(\theta_1)_2=-1,\qquad (\theta_2)_3=-2,\qquad (\theta_3)_4=2.
\end{equation}
Using the results in \eqref{HC2aFa}, \eqref{HC2a-titi} and the classification conditions
\eqref{cl-H}, \eqref{cl-N}, we obtain
\begin{prop}\label{prop-HC2}
The hypercomplex manifold with Hermitian-Norden metrics on a 4-dimensional Lie algebra,
determined by \eqref{HC2a}, belongs to the largest class of the considered manifolds,
i.e. $\HC$, as well as this manifold does not belong to neither $\W_1$ nor $\W_2$ for $J_2$ and $J_3$.
\end{prop}

 The other possibility is the signature of $g$ on $\R$ to be $(-)$, e.g. $e_3\in\R$, where $g(e_3,e_3)=-1$.
By similar computations we establish the same class in the statement of \propref{prop-HC2}.

\subsection{Hypercomplex HN-metric structure of type (hc3)}


We analyze separately the cases of signature (1,1), (0,2) and (2,0) of $g$ on $\g'$.

\subsubsection{}\label{2.3.1.}

Firstly, we consider $g$ of signature (1,1) on $\g'$.

Let us determine $\g$ by %
\begin{equation}\label{HC3a}
[e_2,e_3]=[e_1,e_4]=e_2,\qquad
[e_2,e_1]=[e_4,e_3]=e_4.
\end{equation}
Then we compute covariant derivatives
and the nontrivial ones are %
\begin{equation}\label{HC3a-n}
\begin{array}{ll}
\n_{e_2} e_1=\n_{e_4} e_3=e_4,\qquad &
\n_{e_2} e_2=-\n_{e_4} e_4=e_3,\\[4pt]
\n_{e_2} e_3=-\n_{e_4} e_1=e_2,\qquad & \n_{e_2} e_4=\n_{e_4}
e_2=e_1.
\end{array}
\end{equation}

By virtue of \eqref{HC3a}, \eqref{JJJ} and  \eqref{F}, we obtain that $F_1=0$
and the other components $(F_{\al})_{ijk}$, $\al=2,3$, are as follows
\begin{equation}\label{HC3aFa}
\begin{aligned}
(F_{2})_{212}&=(F_{2})_{221}=(F_{2})_{234}=(F_{2})_{243}=
\\[4pt]
&=-(F_{2})_{414}=(F_{2})_{423}=(F_{2})_{432}=-(F_{2})_{441}=2;
\\[4pt]
(F_{3})_{211}&=-(F_{3})_{222}=-(F_{3})_{233}=(F_{3})_{244}=
\\[4pt]
&=(F_{3})_{413}=(F_{3})_{424}=(F_{3})_{431}=(F_{3})_{442}=-2.
\end{aligned}
\end{equation}
 The only non-zero components of the corresponding Lee forms are
\begin{equation}\label{HC3a-titi}
(\theta_2)_1=(\theta_3)_2=4.
\end{equation}
Using that $F_1=0$, the results in \eqref{HC3aFa}, \eqref{HC3a-titi}
and the classification conditions \eqref{cl-H}, \eqref{cl-N}, we obtain

\begin{prop}\label{prop-HC3}
The hypercomplex manifold with Hermitian-Norden metrics on a 4-dimensional Lie algebra,
determined by \eqref{HC3a}, belongs to the subclass of the K\"ahler manifold
with respect to $J_1$ of the largest class of the considered manifolds, i.e.
\[
\W_0(J_1)\cap\left(\W_1\oplus\W_2\right)(J_2)\cap\left(\W_1\oplus\W_2\right)(J_3),
\]
as well as this manifold does not belong to neither $\W_1$ nor $\W_2$ for $J_2$ and $J_3$.
\end{prop}

\subsubsection{}\label{2.3.2.}

Secondly, we consider $g$ of signature (2,0)  on $\g'$. The case for signature (0,2) is similar.

Let us determine $\g$ by \begin{equation}\label{HC3b}
[e_1,e_3]=[e_4,e_2]=e_1,\qquad [e_1,e_4]=[e_2,e_3]=e_2.
\end{equation} Then we compute covariant derivatives and the
nontrivial ones are \begin{equation}\label{HC3b-n}
\begin{array}{c}
\n_{e_1} e_1=\n_{e_2} e_2=e_3,\qquad
\n_{e_2} e_3=-\n_{e_4} e_1=e_2,\\[4pt]
\n_{e_1} e_3=\n_{e_4} e_2=e_1.
\end{array}
\end{equation}

By virtue of \eqref{HC3b}, \eqref{JJJ} and  \eqref{F}, we obtain the following components of  $(F_{\al})$:
\begin{equation}\label{HC3bFa}
\begin{aligned}
(F_{1})_{114}&=-(F_{1})_{123}=(F_{1})_{132}=-(F_{1})_{141}=
\\[4pt]
&=(F_{1})_{213}=(F_{1})_{224}=-(F_{1})_{231}=-(F_{1})_{242}=-1;
\\[4pt]
(F_{2})_{111}&=(F_{2})_{133}=2,\\[4pt]
(F_{2})_{212}&=(F_{2})_{221}=(F_{2})_{234}=(F_{2})_{243}
\\[4pt]
&=-(F_{2})_{414}=(F_{2})_{423}=(F_{2})_{432}=-(F_{2})_{441}=1;
\\[4pt]
(F_{3})_{222}&=(F_{3})_{233}=2,\\[4pt]
-(F_{3})_{112}&=-(F_{3})_{121}=(F_{3})_{134}=(F_{3})_{143}
\\[4pt]
&=(F_{3})_{413}=(F_{3})_{424}=(F_{3})_{431}=(F_{3})_{442}=-1.
\end{aligned}
\end{equation}
 The only non-zero components of the corresponding Lee forms are
\begin{equation}\label{HC3b-titi}
(\theta_1)_4=-2,\qquad (\theta_2)_1=(\theta_3)_2=4.
\end{equation}
Using the results in \eqref{HC3bFa},  \eqref{HC3b-titi}
and the classification conditions \eqref{cl-H}, \eqref{cl-N},
we obtain that the considered manifold belongs to the class $\W_4(J_1)\cap\W_1(J_2)\cap\W_1(J_3)$.
Remark that, according to \cite{GrMaDi}, necessary and sufficient
conditions a 4-dimensional almost hypercomplex HN-manifold to be
in the class $\W=\W_4(J_1)\cap\W_1(J_2)\cap\W_1(J_3)$ are:
\begin{equation}\label{titaJ}
\theta_2\circ J_2 =\theta_3\circ J_3 =-2\left(\theta_1\circ
J_1\right).
\end{equation}
These conditions are satisfied bearing in mind \eqref{HC3b-titi}.

Let us consider the class $\W^0=\left\{\W\;|\;\mathrm{d}\left(\theta_1\circ J_1\right)=0\right\}$,
which is the class of the (locally) conformally equivalent
$\K$-manifolds, where a conformal transformation of the metric is given by $\bar{g} = e^{2u}g$
for a differential function $u$ on the manifold.

Using \eqref{HC3b-titi} and \eqref{titaJ}, we establish that the considered manifold
belongs to the subclass $\W^0$.

\begin{prop}\label{prop-HC3b}
The hypercomplex manifold with Hermitian-Norden metrics on a 4-dimensional Lie algebra,
determined by \eqref{HC3b}, belongs to the class $\W^0$ of the (locally) conformally
equivalent $\K$-manifolds.
\end{prop}

\subsection{Hypercomplex HN-metric structure of type (hc4)}
\label{HC4}

In this case,   $\g$ is solvable and the derived Lie algebra $\g'$ is 3-dimensional and Abelian.

\subsubsection{}\label{2.4.1.}

Firstly, we fix $e_1\in\g$, for which $g(e_1,e_1)=1$, as an
element ortho\-gonal to $\g'$ with respect to $g$. Therefore $\g$
is determined by \begin{equation}\label{HC4a} [e_1,e_2]=e_2,\qquad
[e_1,e_3]=e_3,\qquad [e_1,e_4]=e_4. \end{equation} Then we compute
covariant derivatives and the nontrivial ones are
\begin{equation}\label{HC4a-n}
\begin{array}{c}
\n_{e_2} e_1=-e_2, \qquad \n_{e_3} e_1=-e_3, \qquad \n_{e_4} e_1=-e_4, \\[4pt]
\n_{e_2} e_2=-\n_{e_3} e_3=-\n_{e_4} e_4=e_1.
\end{array}
\end{equation}

By similar computation as in the previous cases, the components
$(F_{\al})_{ijk}$, $\al=1,2,3$, are as follows:
\begin{equation}\label{HC4aFa}
\begin{aligned}
(F_{1})_{314}&=-(F_{1})_{323}=(F_{1})_{332}=-(F_{1})_{341}=
\\[4pt]
&=-(F_{1})_{413}=-(F_{1})_{424}=(F_{1})_{431}=(F_{1})_{442}=1;
\\[4pt]
(F_{2})_{311}&=(F_{2})_{333}=-2,\\[4pt]
(F_{2})_{214}&=-(F_{2})_{223}=-(F_{2})_{232}=(F_{2})_{241}
\\[4pt]
&=(F_{2})_{412}=(F_{2})_{421}=(F_{2})_{434}=(F_{2})_{443}=-1;
\\[4pt]
(F_{3})_{411}&=(F_{3})_{444}=2,\\[4pt]
(F_{3})_{213}&=(F_{3})_{224}=(F_{3})_{231}=(F_{3})_{242}
\\[4pt]
&=(F_{3})_{312}=(F_{3})_{321}=-(F_{3})_{334}=-(F_{3})_{343}=-1.
\end{aligned}
\end{equation}
 The only non-zero components of the corresponding Lee forms are
\begin{equation}\label{HC4a-titi}
(\theta_1)_2=-(\theta_2)_3=(\theta_3)_4=-2.
\end{equation}
The results in \eqref{HC4aFa}, \eqref{HC4a-titi} and the classification conditions
\eqref{cl-H}, \eqref{cl-N} imply

\begin{prop}\label{prop-HC4}
The hypercomplex manifold with Hermitian-Norden metrics on a 4-dimensional Lie algebra,
determined by \eqref{HC4a}, belongs to the largest class of the considered manifolds,
i.e. $\HC$, as well as this manifold does not belong to neither $\W_1$ nor $\W_2$ for $J_2$ and $J_3$.
\end{prop}

\subsubsection{}\label{2.4.2.}

Secondly, we choose $e_4\in\g$, for which $g(e_4,e_4)=-1$, as an
element orthogonal to $\g'$ with respect to $g$. Therefore, in
this case $\g$ is determined by \begin{equation}\label{HC4b}
[e_4,e_1]=e_1,\qquad [e_4,e_2]=e_2,\qquad [e_4,e_3]=e_3.
\end{equation} Therefore, the nontrivial covariant derivatives are
\begin{equation}\label{HC4b-n}
\begin{array}{c}
\n_{e_1} e_1=\n_{e_2} e_2=-\n_{e_3} e_3=-e_4, \\[4pt]
\n_{e_1} e_4=-e_1, \qquad \n_{e_2} e_4=-e_2, \qquad \n_{e_3}
e_4=-e_3.
\end{array}
\end{equation}

In a similar way we obtain:
\begin{equation}\label{HC4bFa}
\begin{aligned}
(F_{1})_{113}&=(F_{1})_{124}=-(F_{1})_{131}=-(F_{1})_{142}=
\\[4pt]
&=-(F_{1})_{214}=(F_{1})_{223}=-(F_{1})_{232}=(F_{1})_{241}=-1;
\\[4pt]
(F_{2})_{222}&=(F_{2})_{244}=-2,\\[4pt]
(F_{2})_{112}&=(F_{2})_{121}=(F_{2})_{134}=(F_{2})_{143}
\\[4pt]
&=(F_{2})_{314}=-(F_{2})_{323}=-(F_{2})_{332}=(F_{2})_{341}=-1;
\\[4pt]
(F_{3})_{111}&=(F_{3})_{144}=2,\\[4pt]
-(F_{3})_{212}&=-(F_{3})_{221}=(F_{3})_{234}=(F_{3})_{243}
\\[4pt]
&=(F_{3})_{313}=(F_{3})_{324}=(F_{3})_{331}=(F_{3})_{342}=-1.
\end{aligned}
\end{equation}
 The only non-zero components of the corresponding Lee forms are
\begin{equation}\label{HC4b-titi}
(\theta_1)_3=-2,\qquad (\theta_2)_2=-(\theta_3)_1=-4.
\end{equation}

Then, analogously of Case \ref{2.3.2.}, we obtain the following

\begin{prop}\label{prop-HC4b}
The hypercomplex manifold with Hermitian-Norden metrics on a 4-dimensional Lie algebra,
determined by \eqref{HC4b}, belongs to the class $\W^0$ of the (locally) conformally
equivalent $\K$-manifolds.
\end{prop}

\subsection{Hypercomplex HN-metric structure of type (hc5)}
\label{HC5}

In this case,   $\g$ is solvable and $\g'$ is a 3-dimensional Heisenberg algebra.

\subsubsection{}\label{2.5.1.}

Firstly, we fix $e_1\in\g$, for which $g(e_1,e_1)=1$, as an
element ortho\-go\-nal to $\g'$ with respect to $g$. Then $\g$ is
determined by \begin{equation}\label{HC5a} [e_1,e_2]=e_2,\quad
[e_1,e_3]=\frac{1}{2}e_3,\quad [e_1,e_4]=\frac{1}{2}e_4,\quad
[e_3,e_4]=\frac{1}{2}e_2. \end{equation} Then we compute covariant
derivatives and the nontrivial ones are
\begin{equation}\label{HC5a-n}
\begin{aligned}
\n_{e_2} e_2=-2\n_{e_3} e_3=-2\n_{e_4} e_4=e_1,& \\[4pt]
-\n_{e_2} e_1=4\n_{e_3} e_4=-4\n_{e_4} e_3=e_2,& \\[4pt]
-4\n_{e_2} e_4=-2\n_{e_3} e_1=-4\n_{e_4} e_2=e_3,& \\[4pt]
4\n_{e_2} e_3=4\n_{e_3} e_2=-2\n_{e_4} e_1=e_4.&
\end{aligned}
\end{equation}

Analogously of the previous cases we obtain the non-zero
components $(F_{\al})_{ijk}$, $\al=1,2,3$, as follows:
\begin{equation}\label{HC5aFa}
\begin{aligned}
(F_{1})_{314}&=-(F_{1})_{323}=(F_{1})_{332}=-(F_{1})_{341}=
\\[4pt]
&=-(F_{1})_{413}=-(F_{1})_{424}=(F_{1})_{431}=(F_{1})_{442}=\frac{1}{4};
\\[4pt]
(F_{2})_{214}&=-(F_{2})_{223}=-(F_{2})_{232}=(F_{2})_{241}=-\frac{5}{4},
\\[4pt]
(F_{2})_{311}&=-2(F_{2})_{322}=(F_{2})_{333}=-2(F_{2})_{344}=-1,\\[4pt]
(F_{2})_{412}&=(F_{2})_{421}=(F_{2})_{434}=(F_{2})_{443}=-\frac{3}{4};
\\[4pt]
(F_{3})_{213}&=(F_{3})_{224}=(F_{3})_{231}=(F_{3})_{242}=-\frac{5}{4},
\\[4pt]
(F_{3})_{312}&=(F_{3})_{321}=-(F_{3})_{334}=-(F_{3})_{343}=-\frac{3}{4},
\\[4pt]
(F_{3})_{411}&=-2(F_{3})_{422}=-2(F_{3})_{433}=(F_{3})_{444}=1.
\end{aligned}
\end{equation}
 The only non-zero components of the corresponding Lee forms are
\begin{equation}\label{HC5a-titi}
(\theta_1)_2=-\frac{1}{2},\qquad
(\theta_2)_3=-(\theta_3)_4=3.
\end{equation}

The results in \eqref{HC5aFa}, \eqref{HC5a-titi} and the classification conditions
\eqref{cl-H}, \eqref{cl-N} imply

\begin{prop}\label{prop-HC5a}
The hypercomplex manifold with Hermitian-Norden metrics on a 4-dimensional Lie algebra,
determined by \eqref{HC5a}, belongs to the largest class of the considered manifolds,
i.e. $\HC$, as well as this manifold does not belong to neither $\W_1$ nor $\W_2$ for $J_2$ and $J_3$.
\end{prop}

\subsubsection{}\label{2.5.2.}

The other possibility is to choose $e_4\in\g$, for which
$g(e_4,e_4)=-1$, as an element orthogonal to $\g'$ with respect to
$g$. We rearrange the basis in \eqref{HC5a} and then $\g$ is
determined by
\begin{equation}\label{HC5b}
[e_1,e_2]=-\frac{1}{2}e_3,\quad [e_1,e_4]=-\frac{1}{2}e_1, \quad
[e_2,e_4]=-\frac{1}{2}e_2,\quad [e_3,e_4]=-e_3.
\end{equation}

By similar computations we establish the same statement as of \propref{prop-HC5a}
for the Heisenberg algebra introduced by \eqref{HC5b}.


\end{document}